\documentclass[12pt]{amsart}

\usepackage[margin=1in]{geometry}

\usepackage[normalem]{ulem}
\usepackage{amsmath,amscd,amssymb,amsfonts,latexsym,wasysym, mathrsfs, mathtools,hhline,color}
\usepackage[all, cmtip]{xy}
\usepackage{csquotes}
\usepackage{url}
\usepackage{comment}

\definecolor{hot}{RGB}{65,105,225}

\usepackage[pagebackref=true,colorlinks=true, linkcolor=hot ,  citecolor=hot, urlcolor=hot]{hyperref}
\usepackage{ textcomp }
\usepackage{ tipa }
\usepackage{graphicx,enumerate}

\usepackage{geometry}
\theoremstyle{plain}
\newtheorem{theorem}{Theorem}[section]
\newtheorem{prop}[theorem]{Proposition}

\newtheorem{lm}[theorem]{Lemma}

\newtheorem{cor}[theorem]{Corollary}

\newtheorem{lemma}[theorem]{Lemma}
\newtheorem{thrm}[theorem]{Theorem}

\theoremstyle{definition}

\newtheorem{defn}[theorem]{Definition}

\newtheorem{rmk}[theorem]{Remark}

\newtheorem{ex}[theorem]{Example}
\newtheorem*{ex*}{Example}

\def\be{\begin{equation}}
\def\ee{\end{equation}}

\def\bt{\begin{thrm}}
\def\et{\end{thrm}}

\def\bc{\begin{cor}}
\def\ec{\end{cor}}

\def\br{\begin{rmk}}
\def\er{\end{rmk}}

\def\bp{\begin{prop}}
\def\ep{\end{prop}}

\def\bl{\begin{lm}}
\def\el{\end{lm}}

\def\bex{\begin{ex}}
\def\eex{\end{ex}}

\def\bd{\begin{defn}}
\def\ed{\end{defn}}

\newcommand{\CP}{\mathbb{CP}}

\newcommand{\C}{\mathbb{C}}
\newcommand{\Z}{\mathbb{Z}}

\newcommand{\Tor}{\mathrm{Tor}}

\newcommand{\sH}{\mathcal{H}}

\newcommand{\sA}{\mathcal{A}}

\newcommand{\sL}{\mathcal{L}}
\newcommand{\sK}{\mathcal{K}}
\newcommand{\sP}{\mathcal{P}}

\newcommand{\coker}{\mathrm{coker}}

\newcommand{\id}{\mathrm{Id}}

\title[]{Cohomology of $\Z$-local systems on complex hyperplane arrangement complements}

\author{Yongqiang Liu}
\address{Y. Liu: The Institute of Geometry and Physics, University of Science and Technology of China, 96 Jinzhai Road, Hefei 230026 P.R. China} 
\email{liuyq@ustc.edu.cn}
\author{Laurentiu Maxim}
\address{L. Maxim: Department of Mathematics,  University of Wisconsin-Madison,  480 Lincoln Drive, Madison WI 53706-1388, USA,  \newline
{\text and} Institute of Mathematics of the Romanian Academy, P.O. Box 1-764, 70700 Bucharest, ROMANIA.}
\email {maxim@math.wisc.edu}
\author{Botong Wang}
\address{B. Wang: Department of Mathematics,         University of Wisconsin-Madison,  480 Lincoln Drive, Madison WI 53706-1388, USA.}
\email {wang@math.wisc.edu}

\date{\today}

\keywords{hyperplane arrangement, $\Z$-local system, CDO-condition}

\subjclass[2010]{32S22, 52C35, 55N25.}

\begin{document}

\maketitle

\begin{abstract}
We prove a Cohen-Dimca-Orlik type theorem for rank one $\Z$-local systems on complex hyperplane arrangement complements. 
This settles a recent conjecture of S. Sugawara.
\end{abstract}

\section{Introduction}
A finite collection $\sA$ of hyperplanes in $\C^n$ (or $\CP^n$) is called a complex affine (resp., projective) hyperplane arrangement. A fundamental problem in the theory of hyperplane arrangements is to decide whether various invariants of the complement of $\sA$ are determined by the combinatorial information encoded in the intersection poset of the arrangement. For instance, Betti numbers and the cohomology ring of arrangement complements are combinatorially determined (e.g., see \cite{PT}). However, it is still an open question whether the monodromy or the Betti numbers of the Milnor fiber of a central hyperplane arrangement are combinatorially determined; see \cite{PS17} for recent progress in this direction, and also \cite{Dim17} for an overview of the theory. In this note, we show that the cohomology groups of a complex hyperplane arrangement with coefficients in a rank one $\Z$-local system satisfying a certain genericity condition, are combinatorially determined.

The cohomology of local systems on the complement of a complex hyperplane arrangement has been the subject of considerable interest in recent years. For instance, the cohomology of $\C$-local systems is already useful for understanding the Betti numbers of the Milnor fiber associated to a central arrangement \cite{CS95}. Vanishing (or nonresonance) conditions for the cohomology of such $\C$-local systems have been obtained by several authors, including Aomoto \cite{A75}, Esnault-Schectman-Viehweg \cite{ESV}, Kohno \cite{Ko}, Schechtman-Terao-Varchenko \cite{STV}, and Cohen-Dimca-Orlik \cite{CDO03}, see also \cite[Section 6.4]{Dim04} for a survey. More recently, Sugawara  proved a Cohen-Dimca-Orlik type theorem \cite[Theorem 1.3]{Sug22} for the cohomology of a rank one $\Z$-local system on a complexified real hyperplane arrangement. As a consequence of his result, it follows that the cohomology of such local systems are combinatorially determined. His proof uses the chamber structure of the real arrangement complement, and hence it  cannot be extended to complex arrangement complements. However, Sugawara conjectured that the same result should hold for any complex hyperplane arrangement \cite[Conjecture 1.4]{Sug22}.  In this paper, we confirm this conjecture.

\subsection{Setup}
Let $\sA=\{H_1,\cdots,H_d\}$ be a hyperplane arrangement in $\C^n$ with complement \[M=\C^n \setminus \bigcup_{k=1}^d H_k.\] We identify $\C^n$ with $\CP^n \setminus \overline{H}_\infty$, where $ \overline{H}_\infty \subset \CP^n$ denotes the hyperplane at infinity, and let $ \overline{H}_1, \ldots,  \overline{H}_d$ be the projective closures in $\CP^n$ of $H_1, \ldots, H_d$.  Consider the projective arrangement 
\[ \overline{\sA}=\{  \overline{H}_1,\cdots,  \overline{H}_d, \overline{H}_\infty\} \subset \CP^n \] 
and note that $M$ is also the complement of $\overline{\sA}$ in $\CP^n$. 

An edge (or flat) in either $\sA$ or $\overline{\sA}$ is a nonempty intersection of hyperplanes. An edge is said to be dense if the subarrangement of hyperplanes containing it is indecomposable, i.e., these hyperplanes cannot be partitioned into nonempty sets so that, after a change of coordinates, hyperplanes in different sets are in different coordinates. Since this is a combinatorial condition which can be checked in a neighborhood of a given edge, the notion of dense edge can be defined for both affine and projective arrangements.

Let $\sL$ be a rank one $\Z$-local system on $M$, with monodromy representation 
\[
\rho: \pi_1(M)\to \mathrm{Aut}(\Z)\cong\Z^\times.
\]
Since $\Z^\times=\{\pm 1\}$ is abelian, $\rho$ factors through the first homology group $H_1(M,\Z)$ of $M$. 
Therefore, $\sL$ is determined by a $d$-tuple $(t_1,\cdots,t_d)\in \{\pm 1\}^d$, where $t_k=\rho(\gamma_k)$ is the local monodromy of the meridian $\gamma_k$ around the hyperplane $H_k$, $k=1,\cdots,d$. 
Let $t_\infty$ denote the monodromy of the meridian $\gamma_\infty$ around the hyperplane at infinity $\overline{H}_{\infty}$. Since $[\gamma_1]+\cdots+[\gamma_d]+[\gamma_\infty]=0 $ in $H_1(M,\Z)$, it follows that
$t_\infty= (\prod_{k=1}^d t_k)^{-1}$. For any edge $X$ in $\overline{\sA}$, we let  $t_X=\prod_{\overline{H}_k \supseteq X} t_k$, where  the index $k$ is also allowed to be $\infty$.

The following definition concerning rank one $\Z$-local systems is the integral counterpart of the resonance condition (termed ``CDO-condition'' in \cite{Sug22}) introduced by Cohen-Dimca-Orlik \cite{CDO03}.

\bd\label{cdoc} Let $\sL$ be a rank one $\Z$-local system on $M$. We say that $\sL$ satisfies the \emph{CDO-condition} (along $\overline{H}_\infty$) if $t_X\neq 1$ (i.e., $t_X=-1$) for any dense edge $X\subseteq  \overline{H}_\infty$.
\ed 

\subsection{Main result}
The main result of this paper is the following combinatorial formula for the cohomology groups of a rank one $\Z$-local system $\sL$ as in Definition \ref{cdoc}, which was conjectured by Sugawara in \cite[Conjecture 1.4]{Sug22}.

\bt \label{main}
If a rank one $\Z$-local system $\sL$ on $M$ satisfies the CDO-condition, then
\be\label{fm} H^i(M,\sL)\cong
\begin{cases}
\Z^{\beta_n(M)} \oplus \Z_2^{\beta_{n-1}(M)} &  i=n,\\
\Z_2^{\beta_{i-1}(M)} &  1\leq i \leq n-1,\\
0 &\text{ otherwise}.
\end{cases}
\ee
where $\beta_i(M)=  \sum_{j=0}^i (-1)^{i-j} b_j(M)$.
In particular, all the nonzero torsion elements in $H^\bullet(M,\sL)$ are 2-torsion, and the local system cohomology $H^\bullet(M,\sL)$ is combinatorially determined.
\et

Let us mention here that there are examples \cite{ISY22,Sug22} where $\Z/4\Z$-summands appear in $H_1(X,\sL)$ if the $\Z$-local system $\sL$ does not satisfy the CDO-condition. Furthermore, after  tensoring with $\C$, Theorem \ref{main} reduces to the vanishing result of \cite{CDO03} for rank one $\C$-local systems.

\smallskip

Our proof of Theorem \ref{main} makes use of the theory of perverse sheaves (e.g., see \cite{MS, Sch03}), coupled with induction on dimension. We discuss an outline of our strategy of proof in the following subsection.

\subsection{Outline of proof}\label{sec2}
We introduce a new hyperplane $H_0\subset \C^n$ in {\it general position}, and let $\overline{H}_0$ be its closure in $\CP^n$. Denote 
\begin{center} $M':=M\setminus H_0$  \ and \ $M'':=M\cap H_0$.\end{center}

Let $l_0$ and $l_\infty$ be the homogeneous linear functions on $\C^{n+1}$ which define $\overline{H}_0$ and $\overline{H}_\infty$, respectively. By assumption, $l_0$  is chosen generically. Then 
\[
f\coloneqq l_0/l_\infty: M'\to \C^*
\]
is (the restriction of) an affine linear map. Let $\Delta_\infty^*\subset \C^*\subset \CP^1$ be a small punctured disc centered at $\infty$. Let $$U=f^{-1}(\Delta_\infty^*),$$ and let $f_U: U\to \Delta^*_\infty$ be the restriction of $f$ to $U$. In particular, $U$ is open in $M'$. Since the radius of $\Delta_\infty^*$ is small, $f_U$ is a fiber bundle. We denote the fiber of $f_U$ by $M_f$. 

As an intermediate step towards proving Theorem \ref{main}, we show the following result. 
\begin{theorem}\label{thm_2}
If a rank one $\Z$-local system $\sL$ on $M$ satisfies the CDO-condition, then for any $i$, $H^i(U, \sL|_U)$ is  either $0$ or a finite direct sum of $\Z_2$'s.
\end{theorem}

We will use induction on dimension to simultaneously prove Theorems \ref{main} and \ref{thm_2}. 
We apply the inductive hypothesis to the fiber $M_f$ of $f_U$, which is itself an affine hyperplane arrangement complement, after noting that the restriction $\sL|_{M_f}$ of $\sL$ to ${M_f}$ also satisfies the CDO-condition (along the hyperplane at infinity), see Lemma \ref{lemma1} below.

The inductive proof of Theorems \ref{main} and \ref{thm_2} consists of the following steps:
\begin{enumerate}[(i)]
\item\label{item1} Theorem \ref{main} for $\sL|_{M_f}$ on $M_f$ implies Theorem \ref{thm_2} for $\sL|_U$ on $U$. 
\item\label{item2} Theorem \ref{thm_2} for $\sL|_U$ on $U$ implies Theorem \ref{main} for $\sL|_{M'}$ on $M'$. 
\item\label{item3} Theorem \ref{main} for both $\sL|_{M'}$ on $M'$ and $\sL|_{M''}$ on $M''$ implies Theorem \ref{main} as stated (i.e., for $\sL$ on $M$). 
\end{enumerate}

Since dimensions of both $M_f$ and $M''$ are strictly smaller than the dimension of $M$, it follows by induction on dimension that the above statements (i)-(iii) are sufficient to show that both Theorems \ref{main} and \ref{thm_2} hold in all dimensions. 

\smallskip

Our approach in this paper is quite different from that of more classical works \cite{Li, CDO03, DL, Max} etc., where one uses Artin's vanishing theorem and Mayer-Vietoris (or hypercohomology) spectral sequences to infer information about the cohomology of a local system (or of a perverse sheaf) from its
behavior in a small tubular neighborhood of the hyperplane at infinity.  
In order to show that all nonzero torsion elements of $H^\bullet(M,\sL)$ are 2-torsion, we have to distinguish between trivial and non-trivial extensions of $\Z_2$'s, hence Mayer-Vietoris or hypercohomology spectral sequences cannot be used for this purpose. To get a better understanding of the cohomology groups of $\sL$ ``at infinity", we study the neighborhood $U$ defined above, and use its fiber bundle structure over a small puncture disc to control the cohomology groups of $\sL$ on $U$. This is different from using a small tubular neighborhood of the hyperplane at infinity, with its circle-bundle structure, as in the above-mentioned references. Note also that our inductive strategy relies in an essential way on the fiber $M_f$ of the bundle $U$, which is itself the complement of an affine hyperplane arrangement satisfying the CDO-condition.

\section{Monodromy eigenvalues at infinity}
To show Step \eqref{item1} in our proof outline, we need to study the eigenvalues of the monodromy action on  $H^\bullet(M_f, \sL|_{M_f})$. By an eigenvalue of an action on a $\Z$-module we always mean the eigenvalue for the corresponding $\C$-module. So eigenvalues are always complex numbers. 

We begin this section with the following lemma needed for our inductive procedure. In the notations of the previous section, we have the following.
\begin{lemma}\label{lemma1}
The fiber $M_f$ is the complement of an affine hyperplane arrangement. Moreover, the restriction $\sL|_{M_f}$ also satisfies the CDO-condition. 
\end{lemma}
\begin{proof}
By definition, $M_f=M\cap \{l_\infty=\delta l_0\}\subset \CP^n$, where $\delta\in \C$ is nonzero and with sufficiently small absolute value. Since $l_0$ is a general linear form, the hyperplane $\{l_\infty=\delta l_0\}$ is in general position with respect to the hyperplane arrangement $\overline{\sA}=\{ \overline{H}_1,\cdots,  \overline{H}_d, \overline{H}_\infty\}$. Therefore, the restriction $\sL|_{M_f}$ also satisfies the CDO-condition (along the hyperplane at infinity). 
\end{proof}

For Step \eqref{item1} in our proof outline, we need to prove Theorem \ref{thm_2} assuming that Theorem \ref{main} holds for $\sL|_{M_f}$ on $M_f$. We first show that all eigenvalues of the monodromy action on $H^\bullet(M_f, \sL|_{M_f})$ are $-1$, and then deduce Theorem \ref{thm_2}. 

We begin with some lemmas about central arrangements. Let $\mathcal{B}=\{L_0, \ldots, L_m\}$ be a central hyperplane arrangement in $\C^n$, and let $V:=\C^n\setminus (L_0\cup \cdots \cup L_m)$. Let $\sK$ be a rank one local system on $V$, and denote its monodromy around each $L_k$ by $r_k$. Suppose $g_k: \C^n\to \C$ is the linear function whose zero locus is $L_k$. Clearly, $g_0|_V: V\to \C^*$ is a fiber bundle. 
\begin{lemma}\label{lem_eigen}
For any $i$, $R^i(g_0|_V)_{*}(\sK)$ is a local system on $\C^*$ whose monodromy action is given by multiplication by  $r_+=\prod_{k=0}^m r_k$. 
\end{lemma}
\begin{proof}
Fix a loop $\gamma: [0, 1]\to \C^*$ defined by $\gamma(t)=e^{2\pi \sqrt{-1} t}$. Let $V_t:=g_0^{-1}(t)$. Fix a point $\tilde\gamma(0)\in V_{\gamma(0)}$, and let $\tilde\gamma(t)=\gamma(t)\cdot \tilde\gamma(0)$. Then $\tilde\gamma(t)$ is the lifting of the loop $\gamma(t)$ via the map $g_0|_V: V\to \C^*$. Fix a nonzero point $s(0)\in \sK|_{\tilde\gamma(0)}$. Let $s(t)$ be the parallel transport of $s(0)$ along the loop $\tilde\gamma(t)$. Since the element $[\tilde\gamma(t)]\in H_1(V, \Z)$ is equal to the sum of all the meridians around each $L_k$, we have that $s(1)=r_+\cdot s(0)$ in the vector space $\sK|_{\tilde\gamma(0)}$. Therefore, the monodromy action on $H^i(V_{\gamma(0)}, \sK|_{V_{\gamma(0)}})$ induced by going through the loop $\gamma(t)$ is equal to multiplication by $r_+$. 
\end{proof}
\begin{rmk}\label{rmk_eigen}
For a central hyperplane arrangement, the local and global monodromy eigenvalues are the same. Thus, in terms of nearby cycles, the above lemma can be reformulated as follows: if $\mu\colon V\to \C^n$ is the inclusion map, the monodromy action on (any cohomology group of) the nearby cycle complex $\Psi_z Rg_{0*}R\mu_*(\sK)$ is equal to multiplication by $r_+$, where $z$ is the standard coordinate function of $\C$. 
\end{rmk}

\bl \label{lem} Let $\mathcal{B}=\{L_0, \ldots, L_m\}$ be a central hyperplane arrangement in $\C^n$, and let $V$, $\sK$, $r_k$ be defined as in Lemma \ref{lem_eigen}. If the numbers $r_k$ satisfy $\prod_{k=0}^m r_k =-1$, then 
\be\label{v0}
 H^i(V,\sK) \cong H^{i-1}(V^*,\Z_2),  \ \text{for all} \ i,
\ee
where $V^*$ is the projection of $V$ under the Hopf map. 
In particular, each $H^i(V,\sK)$ is either $0$ or a finite direct sum of $\Z_2$'s.
\el 

\begin{proof}
Denote the restriction of the Hopf fibration over $V^*$ by
$$ h\colon V \to V^*,$$
and note that $h$ is a trivial fibration. Since $h$ is a $\C^*$-fiber bundle, $R^j h_* \sK=0$ for $j\neq 0,1$.  
The assumption $\prod_{k=0}^m r_k =-1$ implies that the total turn monodromy operator $T(\sK)$ (cf. \cite[page 210]{Dim04} for a definition) of the local system $\sK$, which in our case is just the monodromy of $\sK$ in the fiber $\C^*$ of $h$,  is the multiplication by $-1$.

Note that $\ker(T(\sK)-\id)=0$ and $\coker(T(\sK)-\id)\cong \Z_2$. As in \cite[Proposition 6.4.3]{Dim04}, we get that $R^0 h_*\sK=0$ and $R^1 h_* \sK$ is a rank one $\Z_2$-local system (hence a constant sheaf, since $\Z_2$ has a trivial automorphism group). 
Using the Leray spectral sequence for $h$, we then have 
$$H^i(V,\sK)\cong H^{i-1}(V^*, R^1 h_* \sK) \cong H^{i-1}(V^*, \Z_2) ,$$
 which implies that $H^i(V,\sK) $ is either $0$ or a finite direct sum of $\Z_2$'s.
\end{proof}

\bc \label{cor central} Let $\mathcal{B}=\{L_0, \ldots, L_m\}$ be a central hyperplane arrangement in $\C^n$ with complement $V$. Suppose $\sK$ is a rank one $\Z$-local system on $V$ such that $t_X=-1$ for any dense edge $X$ of $\sA$ contained in $L_0$. Then $H^\bullet(V,\sK)$ is either $0$ or a finite direct sum of $\Z_2$'s.
\ec
\begin{proof}

Consider the decomposition of the complement $V$ into a finite product of indecomposable central hyperplane arrangement complements, say $V=V_1\times \cdots \times V_r$, with a corresponding rank one $\Z$-local system $\sK_j$ on $V_j$ for $1\leq j\leq r$ such that $$\sK \cong p_1^*\sK_1 \otimes \cdots \otimes p_r^*\sK_r.$$ Here, $p_j$ is the  projection map to the $j$-th factor. Without loss of generality, we can assume that the hyperplane $L_0$ comes from $V_1$. By assumption, the turn monodromy of $\sK$ around the dense edge $0\times \overline{V}_2\times\cdots \times\overline{V}_r$ is equal to $-1$, where each $\overline{V}_i\subset \C^n$ is the closure of $V_i$. Hence, the total turn monodromy of $\sK_1$ on $V_1$ is $-1$. By the previous lemma, $H^\bullet(V_1, \sK_1)$ is either $0$ or a finite direct sum of $\Z_2$'s, and hence, by K\"unneth formula, $H^\bullet(V, \sK)$ is either $0$ or a finite direct sum of $\Z_2$'s. 
\end{proof}

We can now prove the following result about the monodromy eigenvalues of the $\C^*$-action on $H^\bullet(M_f, \sL|_{M_f})$.

\begin{prop}
Under the notations and assumptions of Theorem \ref{thm_2}, the monodromy eigenvalues of $H^\bullet(M_f, \sL|_{M_f})$ are all $-1$. 
\end{prop}

\begin{proof}
Let 
\[
\pi\colon W=\mathrm{Bl}_{\overline{H}_\infty\cap \overline{H}_0}\CP^n\to \CP^n
\]
be the map defined by the blowup of $\CP^n$ along $\overline{H}_\infty\cap \overline{H}_0$.
 
By construction, we have a pencil $\bar{f}\colon W\to \CP^1$, which extends the linear map $f=l_0/l_\infty\colon \C^n\to \C$.
Since $\pi\colon W\to \CP^n$ is an isomorphism over $M'$, the inclusion $M'\subset \CP^n$ lifts to an open embedding $\iota\colon M' \hookrightarrow   W$. 

The stalk of the local system $\sH^i(R\bar{f_*}R\iota_*(\sL|_{M'}))|_{\Delta_\infty^*}$ is isomorphic to $H^i(M_f, \sL|_{M_f})$, so the assertion of the proposition is equivalent to showing that the monodromy eigenvalues of this local system are $-1$ for all $i$. In other words, we need to show that the only monodromy eigenvalue of the nearby cycle complex $\Psi_{1/z}(R\bar{f_*}R\iota_*(\sL|_{M'}))$ is $-1$, where $z$ is the standard coordinate function on $\C=\CP^1\setminus \{\infty\}$ and $1/z$ is a coordinate function on $\CP^1\setminus \{0\}$. Since the nearby cycle functor commutes with proper pushforward, it suffices to show that the monodromy action on the stalks of the nearby cycle complex $\Psi_{1/\bar{f}}(R\iota_*(\sL|_{M'}))$ has only eigenvalue $-1$. Here we regard $1/\bar{f}$ as a $\C$-valued function $1/\bar{f}: W\setminus \bar{f}^{-1}(0)\to \C$. 

Since $\sL$ satisfies the CDO-condition (along $\overline{H}_\infty$), by using Lemma \ref{lem_eigen}  we see that, away from the exceptional divisor $E$, the monodromy action on the stalks of the nearby cycle complex $\Psi_{1/\bar{f}}(R\iota_*(\sL|_{M'}))$ has only the eigenvalue $-1$. 

Along $E\cap \bar{f}^{-1}(\infty)$, we claim that the stalks of the nearby cycle complex $\Psi_{1/\bar{f}}(R\iota_*(\sL|_{M'}))$ consist of torsion elements only, which means that it has no monodromy eigenvalue (in $\C$). In fact, given any point $x\in E\cap \bar{f}^{-1}(\infty)$, let $M_{f, x}$ be the local Milnor fiber of $f\colon M'\to \C^*$ near $x$. More precisely, $M_{f, x}=f^{-1}(1/\delta)\cap B_{\epsilon}(x)$, where $B_{\epsilon}(x)\subset W$ is a ball centered at $x$ with radius $\epsilon$, and $\delta, \epsilon\in \C$ satisfy $0<|\delta|\ll |\epsilon|\ll 1$. By the construction of blowup, the local Milnor fiber can be expressed as
\begin{equation}\label{eq_M}
M_{f, x}\cong M'\cap B_{\epsilon}(\pi(x))\cap \{l_\infty=\delta l_0\},
\end{equation}
where $B_{\epsilon}(\pi(x))\subset \CP^n$ is a ball centered at $\pi(x)$ of radius $\epsilon$, with 
$l_0, l_\infty$ the defining linear forms of $\overline{H}_0, \overline{H}_\infty$, respectively, and $\delta, \epsilon\in \C$ satisfying $0<|\delta|\ll |\epsilon|\ll 1$. The cohomology of the stalk of $\Psi_{1/\bar{f}}(R\iota_*(\sL|_{M'}))$ at $x$ is isomorphic to $H^\bullet(M_{f, x}, \sL|_{M_{f, x}})$. Since $l_0$ is a general linear form and since $\delta\in \C$ is also general, the hyperplane $\overline{H}_{\delta}\coloneqq\{l_\infty=\delta l_0\}$ is in general position. Therefore, $M'\cap \overline{H}_\delta$ is a hyperplane arrangement complement in $\overline{H}_\delta$, and $\sL|_{M'\cap \overline{H}_\delta}$ satisfies the CDO-condition along $\overline{H}_\infty\cap \overline{H}_\delta$. Since $M_{f, x}$ is the intersection of $M'\cap \overline{H}_\delta$ and a small ball in $\overline{H}_\delta$, $M_{f, x}$ is homeomorphic to a central hyperplane arrangement complement.
Then it follows from Corollary \ref{cor central} that $H^\bullet(M_{f, x}, \sL|_{M_{f, x}})$ consists only of torsion elements.
\end{proof}


\section{Cohomology groups at infinity} 
We can now proceed with Step \eqref{item1} in our proof of Theorem \ref{main}.
\begin{proof}[Proof of Theorem \ref{thm_2} assuming Theorem \ref{main} for $M_f$]
We have seen in Lemma \ref{lemma1} that the restriction $\sL|_{M_f}$ of $\sL$ to the fiber $M_f$ satisfies the CDO-condition (along the hyperplane at infinity). Assume that  Theorem \ref{main} holds for $\sL|_{M_f}$ on $M_f$. In particular, the torsion subgroup of $H^\bullet(M_f, \sL|_{M_f})$ consists only of $2$-torsion.

Consider the Wang sequence
\begin{equation}\label{eq_long}
\cdots\to H^{i-1}(M_f, \sL|_{M_f})\xrightarrow[]{T-1} H^{i-1}(M_f, \sL|_{M_f})\to H^{i}(U, \sL|_U)\to H^{i}(M_f, \sL|_{M_f})\to \cdots
\end{equation}
where $T: H^\bullet(M_f, \sL|_{M_f})\to H^\bullet(M_f, \sL|_{M_f})$ is the monodromy action. 

In the above long exact sequence, let us denote the image of $H^{i-1}(M_f, \sL|_{M_f})$ in $H^i(U, \sL|_U)$ by $A^i$, and the image of $H^i(U, \sL|_U)$ in $H^{i}(M_f, \sL|_{M_f})$ by $B^i$. Then, we have a short exact sequence
\begin{equation}\label{eq_short}
0\to A^i\to H^i(U, \sL|_U)\to B^i\to 0. 
\end{equation}
 Notice that $A^i$ is isomorphic to the cokernel of $T-1: H^{i-1}(M_f, \sL|_{M_f})\to H^{i-1}(M_f, \sL|_{M_f})$. Since the only eigenvalue of $T$ is $-1$ and $H^{i-1}(M_f, \sL|_{M_f})$ has only 2-torsion, it follows that 
 \begin{equation}\label{eq_card1}
|A^i|\leq |H^{i-1}(M_f, \sL|_{M_f})\otimes_\Z \Z_2|,
\end{equation}
where $|\cdot |$ denotes the cardinality.
Similarly, $B^i$ is isomorphic to the kernel of $T-1: H^{i}(M_f, \sL|_{M_f})\to H^{i}(M_f, \sL|_{M_f})$, which is contained in $H^{i}(M_f, \sL|_{M_f})_{tors}$. Therefore, 
\begin{equation}\label{eq_card2}
|B^i|\leq |H^{i}(M_f, \sL|_{M_f})_{tors}|. 
\end{equation}
By combining the inequalities \eqref{eq_card1} and \eqref{eq_card2}, and using the short exact sequence \eqref{eq_short}, we obtain
\begin{equation}\label{eq_leq}
|H^i(U, \sL|_U)|\leq |H^{i-1}(M_f, \sL|_{M_f})\otimes_\Z \Z_2|\cdot |H^{i}(M_f, \sL|_{M_f})_{tors}|. 
\end{equation}
In particular, $H^\bullet(U, \sL|_U)$ consists only of torsion elements.

On the other hand, notice that $\sL\otimes_{\Z}\Z_2$ is isomorphic to the trivial rank one $\Z_2$-local system. 
We  claim that the induced monodromy action on $H^1(M_f, \Z_2)$ is trivial. In fact, if we denote by $l_i$ the homogeneous linear function defining $\overline{H}_i$, and consider the composition
\[
M_f\to M'\xrightarrow{l_i/l_\infty} \C^*,
\]
where the first map is the inclusion,  then the pullback of a generator of $H^1(\C^*, \Z_2)$ to $H^1(M_f, \Z_2)$ is fixed by the induced monodromy action. Since $H^1(M_f, \Z_2)$ is generated by such pullbacks, the induced monodromy action on $H^1(M_f, \Z_2)$ is trivial.  
Since the cohomology ring $H^\bullet(M_f, \Z_2)$ is generated in degree one, the monodromy action on $H^\bullet(M_f, \Z_2)$ is also trivial. Now, using the Wang sequence with $\Z_2$ coefficients, we obtain
\begin{equation}\label{eq_UM}
H^i(U, \sL|_U\otimes_\Z \Z_2)\cong H^i(U, \Z_2)\cong H^i(M_f, \Z_2)\oplus H^{i-1}(M_f, \Z_2).
\end{equation}
Moreover, by the universal coefficient theorem, 
\begin{equation*}
H^i(U, \sL|_U\otimes_\Z \Z_2)\cong H^i(U, \sL|_U)\otimes_\Z \Z_2\oplus \Tor_1^\Z(H^{i+1}(U, \sL|_U), \Z_2),
\end{equation*}
and, since we assumed that  $H^\bullet(M_f, \sL|_{M_f})_{tors}$ is all $2$-torsion,
\begin{align*}
H^i(M_f, \Z_2)&\cong H^i(M_f, \sL|_{M_f})\otimes_\Z \Z_2\oplus  \Tor_1^\Z(H^{i+1}(M_f, \sL|_{M_f}), \Z_2)\\
&\cong H^i(M_f, \sL|_{M_f})\otimes_\Z \Z_2\oplus  H^{i+1}(M_f, \sL|_{M_f})_{tors}.
\end{align*}
Thus, the isomorphism \eqref{eq_UM} can be rewritten as
\begin{multline*}\label{eq_iso}
H^i(U, \sL|_U)\otimes_\Z \Z_2\oplus \Tor_1^\Z(H^{i+1}(U, \sL|_U), \Z_2)\\
\cong H^i(M_f, \sL|_{M_f})\otimes_\Z \Z_2 \oplus H^{i-1}(M_f, \sL|_{M_f})\otimes_\Z \Z_2 \oplus  H^{i+1}(M_f, \sL|_{M_f})_{tors}\oplus H^{i}(M_f, \sL|_{M_f})_{tors}. 
\end{multline*}
Notice also that
\begin{equation}\label{eq_2geq}
|H^i(U, \sL|_U)|\geq |H^i(U, \sL|_U)\otimes_\Z \Z_2|  \quad\text{and}\quad |H^{i+1}(U, \sL|_U)|\geq |\Tor_1^\Z(H^{i+1}(U, \sL|_U), \Z_2)|.
\end{equation}
Therefore, 
\begin{multline}\label{eq_geq}
|H^i(U, \sL|_U)|\cdot |H^{i+1}(U, \sL|_U)|\\
\geq |H^{i-1}(M_f, \sL|_{M_f})\otimes_\Z \Z_2| \cdot |H^i(M_f, \sL|_{M_f})\otimes_\Z \Z_2| \cdot |H^{i}(M_f, \sL|_{M_f})_{tors}| \cdot   |H^{i+1}(M_f, \sL|_{M_f})_{tors}|. 
\end{multline}
Comparing \eqref{eq_leq} and \eqref{eq_geq}, we see that both of these inequalities must be equalities. Hence, the inequalities in \eqref{eq_2geq} must also be equalities. In particular, $$|H^i(U, \sL|_U)|= |H^i(U, \sL|_U)\otimes_\Z \Z_2|,$$ which means that $H^i(U, \sL|_U)$ consists only of $2$-torsion elements. 
\end{proof}


\section{Torsion in the cohomology groups}
In this section, we discuss Steps \eqref{item2} and \eqref{item3} of the proof of Theorem \ref{main}, as outlined in Subsection \ref{sec2}. We employ the notations from Subsection \ref{sec2}.

For proving Step \eqref{item2}, it suffices to show the following result. 
\begin{prop}\label{prop_UM}
Assume that $H^k(U, \sL|_U)$ consists of 2-torsions only for every $k$. Then, we have natural isomorphisms 
\[
H^k(M', \sL|_{M'})\cong H^k(U, \sL|_{U})\quad \text{for} \quad k<n
\]
and a natural injective map
\[
H^n(M', \sL|_{M'})_{tors}\hookrightarrow H^n(U, \sL|_{U}).
\]
\end{prop}

Before proving Proposition \ref{prop_UM}, let us recall the following result about generic projections, see, e.g., \cite{LMW}. 

\begin{theorem}\label{thm_constant}
Let $\sP$ be a perverse sheaf on $\C^n$ over a ring $R$. Let $g: \C^n\to \C$ be a general linear function. Then, up to shifts of constant sheaves on $\C$, $Rg_*(\sP)$ is a perverse sheaf. In other words, the perverse cohomology sheaf $^p\sH^k(Rg_*(\sP))$ is a shift of a constant sheaf on $\C$, for all $k\neq 0$. Moreover, if $R$ is a Dedekind domain and $\sP$ is strongly perverse, then $^{p^+}\sH^k(Rg_*(\sP))$ is the shift of a constant sheaf on $\C$, for all $k\neq 0$, where $p^+$ is the dual of the standard $t$-structure $p$. 
\end{theorem}

Recall here, e.g., see \cite[Section 10.2.3]{MS}, that if $\left({^pD}^{\leq 0}(X;R), {^pD}^{\geq 0}(X;R)\right)$ is the perverse $t$-structure on the derived category $D^b_c(X;R)$ of bounded constructible $R$-complexes on a variety $X$, then the dual perverse $t$-structure $p^+$ on $D^b_c(X;R)$ is defined by
$$^{p^+}D^{\leq 0}(X;R):=\mathcal{D}\left({^pD}^{\geq 0}(X;R)\right) \quad \text{and} \quad  
\, ^{p^+}D^{\geq 0}(X;R):=\mathcal{D}\left({^pD}^{\leq 0}(X;R)\right), $$
where $\mathcal{D}$ denotes the Verdier duality functor on $X$. If $R$ is a Dedekind domain, we say that $\sP \in D^b_c(X;R)$ is strongly perverse if $\sP \in {^pD}^{\leq 0}(X;R) \cap \,^{p^+}D^{\geq 0}(X;R)$. We refer to \cite[Proposition 10.2.49, Proposition 10.2.56]{MS} for a characterization of the dual perverse $t$-structure in terms of the perverse $t$-structure when $R$ is a Dedekind domain. For instance, one has in this case the inclusions
$${^{p}D}^{\leq 0}(X;R)\subset \,^{p^+}D^{\leq 0}(X;R)\subset {^{p}D}^{\leq 1}(X;R) \quad \text{and} \quad  ^{p^+}D^{\geq 0}(X;R)\subset {^{p}D}^{\geq 0}(X;R).$$

\begin{proof}[Proof of Proposition \ref{prop_UM}]
Consider the open embeddings $j_0\colon \C^*\to \C_0:=\CP^1\setminus \{\infty\}$ and $j_\infty\colon  \C^*\to \C_\infty:=\CP^1\setminus \{0\}$, together with the closed embedding  $i\colon \{\infty\}\to \C_\infty$. 
Applying the attaching triangle
\be\label{att} {j_\infty}_!j_\infty^* \to id \to Ri_{ *} i^* \to\ee
to the complex $R{j_\infty}_*Rf_*(\sL|_{M'}[n])$, one has the adjunction distinguished triangle:
\begin{equation}\label{eq_tri}
j_{\infty !}Rf_*(\sL|_{M'}[n])\to Rj_{\infty *}Rf_*(\sL|_{M'}[n])\to Ri_{*} i^* Rj_{\infty*}Rf_*(\sL|_{M'}[n])\xrightarrow{[1]}.
\end{equation}
Notice that we have natural isomorphisms 
\begin{equation*}
\begin{split}
H^\bullet(\C_{\infty}, Ri_{*} i^* Rj_{\infty*}Rf_*(\sL|_{M'}[n]))& \cong \mathcal{H}^\bullet(Rj_{\infty*}Rf_*(\sL|_{M'}[n]))_\infty \\
& \cong H^\bullet(\Delta_\infty, Rj_{\infty *}Rf_*(\sL|_{M'}[n])|_{\Delta_\infty}) \\
& \cong H^\bullet(\Delta_\infty^*, Rf_*(\sL|_{M'}[n])|_{\Delta_\infty^*})
\\ &\cong H^\bullet(U, \sL|_{U}[n])
\end{split}
\end{equation*}
and
\[
H^\bullet(\C_\infty, Rj_{\infty*}Rf_*(\sL|_{M'}[n]))
\cong H^\bullet(\C^*, Rf_*\sL|_{M'}[n]) 
\cong H^\bullet(M', \sL|_{M'}[n]). 
\]
Thus, the cohomology long exact sequence associated to \eqref{eq_tri} is of the form
\begin{equation}\label{eq_long2}
\cdots\to H^k(\C_\infty, j_{\infty!}Rf_*(\sL|_{M'}[n]))\to H^k(M', \sL|_{M'}[n])\to H^k(U, \sL|_{U}[n])\to\cdots.
\end{equation}

Note that for any constant sheaf $\mathcal{E}$ of abelian groups on $\C^*$, we have $H^k(\C_\infty, j_{\infty!}(\mathcal{E}))=0$ for all $k$. Indeed, applying the attaching triangle \eqref{att} to the complex 
$Rj_{\infty*}(\mathcal{E})$ on $\C_\infty$, the associated cohomology long exact sequence becomes
$$ \cdots \to H^k(\C_\infty, j_{\infty!}(\mathcal{E})) \to H^k(\C_\infty, Rj_{\infty*}(\mathcal{E})) \to \mathcal{H}^k(Rj_{\infty*}(\mathcal{E}))_\infty \to \cdots$$
with the arrow $H^k(\C_\infty, Rj_{\infty*}(\mathcal{E})) \to \mathcal{H}^k(Rj_{\infty*}(\mathcal{E}))_\infty $ being an isomorphism for all $k$ (by constructibility and homotopy equivalence). 
Also notice that  the constructible complex $Rf_*(\sL|_{M'}[n])$ on $\C^*$ is the restriction of the constructible complex $Rf_*(\sL|_{M}[n])$ to $\C^*$, where we abuse the notations and use the first and second $f$ to denote $f: M'\to \C^*$ and $f: M\to \C_0=\CP^1\setminus \{\infty\}$, respectively (see, e.g., \cite[Prop.10.7(4)]{Bo}). 
Therefore, by the second part of Theorem \ref{thm_constant}, together with the perverse hypercohomology spectral sequence for $p^+$, and the fact that $j_\infty$ is a quasi-finite affine morphism (hence $j_{\infty!}$ preserves the dual perverse $t$-structure, e.g., see \cite[Corollary 10.3.30, Theorem 10.3.69]{MS}), we get 
\begin{equation}\label{eq_iso1}
H^k(\C_\infty, j_{\infty!}Rf_*(\sL|_{M'}[n]))\cong  H^k(\C_\infty, j_{\infty!}\, ^{p^+}\sH^0(Rf_*(\sL|_{M'}[n]))).
\end{equation}

Let us next consider the commutative diagram
$$\xymatrix{
\C^* \ar[r]^{j_{0}} \ar[d]_{j_{\infty}} & \C_0 \ar[d]^{i_0}\\
\C_\infty  \ar[r]_{i_\infty} & \CP^1 
}$$
and note that (cf. \cite[Lemma 6.0.5]{Sch03})
 \begin{equation}\label{basec}
i_{0!} Rj_{0*} \simeq Ri_{\infty *} j_{\infty !}.
\end{equation}
Hence
\begin{equation}\label{eq_iso2}
\begin{split}
H^k(\C_\infty, j_{\infty!}\, ^{p^+}\sH^0(Rf_*(\sL|_{M'}[n]))) &\cong 
H^k(\CP^1, Ri_{\infty *} j_{\infty!}\, ^{p^+}\sH^0(Rf_*(\sL|_{M'}[n]))) \\ &\cong
H^k(\CP^1, i_{0!} Rj_{0*}\, ^{p^+}\sH^0(Rf_*(\sL|_{M'}[n]))) \\ &\cong
H^k_c(\C_0, Rj_{0*}\, ^{p^+}\sH^0(Rf_*(\sL|_{M'}[n]))).
\end{split}
\end{equation}
Since $j_0$ is a quasi-finite affine morphism, we have that $Rj_{0*}\, ^{p^+}\sH^0(Rf_*(\sL|_{M'}[n]))$ is perverse with respect to the dual perverse $t$-structure (e.g., see \cite[Corollary 10.3.30, Theorem 10.3.69]{MS}).
Applying Artin's vanishing theorem (e.g., see \cite[Theorem 10.3.59(2)]{MS}) for $H^k_c(\C_0, Rj_{0*}\, ^{p^+}\sH^0(Rf_*(\sL|_{M'}[n])))$, and by the description of complexes in $^{p^+}D^{\geq 0}_c(pt, \Z)$ in \cite[Proposition 10.2.49]{MS}, we have 
\[
H^k_c(\C_0, Rj_{0*}\, ^{p^+}\sH^0(Rf_*(\sL|_{M'}[n])))=0 \quad\text{for}\quad k<0,
\]
and $H^0_c(\C_0, Rj_{0*}\, ^{p^+}\sH^0(Rf_*(\sL|_{M'}[n])))$ is torsion-free.
By the  isomorphisms \eqref{eq_iso1} and \eqref{eq_iso2}, we have 
\begin{equation}\label{eq_vanish}
H^k(\C_\infty, j_{\infty!}Rf_*(\sL|_{M'}[n]))=0 \quad\text{for}\quad k<0,
\end{equation}
and $H^0(\C_\infty, j_{\infty!}Rf_*(\sL|_{M'}[n]))$ is torsion-free. 
 
By \eqref{eq_vanish}, the long exact sequence \eqref{eq_long2} gives isomorphisms
\[
H^k(M', \sL|_{M'}[n]))\cong H^k(U, \sL|_{U}[n])\quad\text{for}\quad k<-1,
\]
and an exact sequence
\begin{multline}
0\to H^{-1}(M', \sL|_{M'}[n]))\to H^{-1}(U, \sL|_{U}[n])\\
\to H^0(\C_\infty, j_{\infty!}Rf_*(\sL|_{M'}[n]))\to H^{0}(M', \sL|_{M'}[n]))\to H^{0}(U, \sL|_{U}[n]).
\end{multline}
Recall that $H^0(\C_\infty, j_{\infty!}Rf_*(\sL|_{M'}[n]))$ is torsion-free. Moreover, by Theorem~\ref{thm_2}, we have that $H^{-1}(U, \sL|_{U}[n])$ is a torsion group. Thus, the map $$H^{-1}(U, \sL|_{U}[n])
\to H^0(\C_\infty, j_{\infty !}Rf_*(\sL|_{M'}[n]))$$ is the zero map. Hence we have
\[
H^{-1}(M', \sL|_{M'}[n]))\cong H^{-1}(U, \sL|_{U}[n])
\]
and the natural map 
\[
H^0(M', \sL|_{M'}[n])_{tors}\to H^0(U, \sL|_{U}[n])
\]
is injective. After moving the shifts to the cohomology degrees, we have proved the desired statements. 
\end{proof}

Let us finally justify Step \eqref{item3} of the proof of Theorem \ref{main}.
Assuming that Theorem \ref{main} holds for $\sL|_{M''}$ and $\sL|_{M'}$, we need to prove that Theorem 1.2 holds for $\sL$ on ${M}$. We first show that all the torsion elements of $H^\bullet(M, \sL)$ are 2-torsions.

Consider the following ``deletion-restriction'' long exact sequence for $M=M' \cup M''$ (see e.g., \cite[p.222]{Dim04}):
\[
\cdots\xrightarrow[]{\delta^{i-1}} H^{i-2}(M'', \sL|_{M''})\to H^{i}(M, \sL)\to H^{i}(M', \sL|_{M'})\xrightarrow[]{\delta^i} H^{i-1}(M'', \sL|_{M''})\to \cdots.
\]
Since Theorem \ref{main} is assumed to hold for $\sL|_{M'}$, in order to show that all the torsion elements of $H^\bullet(M, \sL)$ are 2-torsions, it suffices to show that the residue maps $\delta^i$ are surjective for $1\leq i\leq n$. This will also show that the only free part in $H^\bullet(M, \sL)$ may appear in degree $n$. 

Since $H^{n+1}(M, \sL)=0$, the above long exact sequence implies that $\delta^n$ is surjective. To show that the other $\delta^i$ are surjective, we consider the following commutative diagram with exact rows, induced by the functoriality of the universal coefficient theorem (using the fact that $\sL\otimes_\Z \Z_2$ is the trivial rank one $\Z_2$-local system):
\[
\xymatrix{
 H^{i-1}(M', \Z_2)\ar[d] \ar[r]&\Tor_1^{\Z}(H^{i}(M', \sL|_{M'}), \Z_2)\ar[d]\ar[r]&0\\
H^{i-2}(M'', \Z_2)\ar[r]& \Tor_1^{\Z}(H^{i-1}(M'', \sL|_{M''}), \Z_2)\ar[r]&0 .
}
\]
 Notice that the left vertical map is the standard residue map with $\Z_2$-coefficients, which is surjective (e.g., as it follows from the proof of Orlik-Solomon theorem). Thus, the right vertical map must also be surjective. Since Theorem \ref{main} is assumed to hold for both $\sL|_{M'}$ and $\sL|_{M''}$, it follows that for $i<n$ both $H^{i}(M', \sL|_{M'})$ and $H^{i-1}(M'', \sL|_{M''})$ consist only of 2-torsion elements. Thus, the surjectivity of 
\[
\Tor_1^{\Z}(H^{i}(M', \sL|_{M'}), \Z_2)\to \Tor_1^{\Z}(H^{i-1}(M'', \sL|_{M''}), \Z_2)
\]
implies the surjectivity of
\[
\delta^i: H^{i}(M', \sL|_{M'})\to H^{i-1}(M'', \sL|_{M''}),
\]
as desired.

\smallskip

To finish the proof of Theorem \ref{main}, 
we use again the fact that $\sL\otimes_\Z \Z_2$ is the trivial rank one $\Z_2$-local system on $M$. Since $H^*(M,\Z)$ is free, we have for any degree $i$ that 
$$\dim_{\Z_2} H^i(M, \sL\otimes_\Z \Z_2)=\dim_{\Z_2} H^i(M, \Z_2) = {\rm rank} H^i(M,\Z) =b_i(M).$$
So the desired formula \eqref{fm} for the cohomology groups of $\sL$ on $M$ follows from the universal coefficient theorem.

\bex Let $\sA$ be a central hyperplane arrangement in $\C^2$ with $d$ hyperplanes. Let $\sL$ be a rank one $\Z$-local system on its complement $M$ and assume that $\prod_{k=1}^d t_i =-1$. 
There is only one dense edge contained in $ \overline{H}_\infty$, which is $ \overline{H}_\infty$ itself. The corresponding monodromy is $-1$, hence $\sL$ satisfies the CDO-condition.
Then we get by Theorem \ref{main} that 
$$H^i(M,\sL)=\begin{cases}
\Z_2^{d-1}, &  i=2,\\
\Z_2, &  i=1,\\
0, &\text{otherwise.}
\end{cases}$$
\eex

\end{document}